\documentclass[12pt]{amsart}
\usepackage{amscd,amssymb}

\usepackage[graph,frame,poly,arc]{xy}  
\usepackage[plainpages,backref,urlcolor=blue]{hyperref}

\theoremstyle{plain}
\newtheorem{thm}[subsection]{Theorem}

\theoremstyle{definition}
\newtheorem{rk}[subsection]{Remark}

\numberwithin{equation}{section}
\newcommand{\M}{{\mathcal M}}

\newcommand{\C}{\mathbb{C}}
\newcommand{\PP}{\mathbb{P}}

\DeclareMathOperator{\rank}{rank}
\DeclareMathOperator{\im}{im}
\DeclareMathOperator{\dd}{d}

\begin{document}

\title [A minimal resolution for the Jacobian ideal]
{A minimal resolution for the Jacobian ideal of a generic curve arrangement}

\author[Alexandru Dimca]{Alexandru Dimca}
\address{Universit\'e C\^ ote d'Azur, CNRS, LJAD, France and Simion Stoilow Institute of Mathematics,
P.O. Box 1-764, RO-014700 Bucharest, Romania}
\email{Alexandru.Dimca@univ-cotedazur.fr}

\author[Gabriel Sticlaru]{Gabriel Sticlaru}
\address{Faculty of Mathematics and Informatics,
Ovidius University
Bd. Mamaia 124, 900527 Constanta,
Romania}
\email{gabriel.sticlaru@gmail.com }

\subjclass[2010]{Primary 14H50; Secondary  14B05, 13D02, 32S22}

\keywords{Jacobian ideal, Jacobian syzygy module, exponents}

\begin{abstract} We consider a nodal curve $C$ in the complex projective plane whose irreducible components $C_i$ are smooth. A minimal set of generators $G$ for the first and second syzygy modules of the Jacobian ideal of $C$  are described, using recent results by  Th. Kahle, H. Schenck,  B. Sturmfels and M. Wiesmann on the likelihood correspondence.
The elements of $G$ have explicit formulas in terms of the equations $f_i=0$ of the irreducible components $C_i$ of $C$. Similar results, including extensions to hypersurfaces arrangements in $\PP^n$ were obtained by R. Burity, Z. Ramos, A. Simis and  St. Toh\u aneanu with a
genericity assumption which may not be easy to test in practice.

\end{abstract}
 
\maketitle

\section{Introduction} 

Let $S=\C[x,y,z]$ be the polynomial ring in three variables $x,y,z$ with complex coefficients, and let $C:f=0$ be a reduced curve of degree $d\geq 3$ in the complex projective plane $\PP^2$. 
We denote by $J_f$ the Jacobian ideal of $f$, i.e. the homogeneous ideal in $S$ spanned by the partial derivatives $f_x,f_y,f_z$ of $f$, and  by $M(f)=S/J_f$ the corresponding graded quotient ring, called the Jacobian (or Milnor) algebra of $f$.
Consider the graded $S$-module of Jacobian syzygies of $f$, namely
$$D_0(C)=D_0(f)=\{(a,b,c) \in S^3 \ : \ af_x+bf_y+cf_z=0\}.$$
We say that $C:f=0$ is an {\it $s$-syzygy curve} if  the module $D_0(f)$ is minimally generated by $s$ homogeneous syzygies, say $r_1,r_2,...,r_s$, of degrees $d_j=\deg r_j$ ordered such that $$1 \leq d_1\leq d_2 \leq ...\leq d_s.$$ 
We call these degrees the {\it exponents} of the curve $C$, and $r_1,...,r_s$ as a {\it minimal set of generators } for the module  $D_0(f)$. It is known that
\begin{equation}
\label{boundG}
d_1\leq d_2 \leq d_3 \leq d-1,
\end{equation}
see for instance \cite[Theorem 2.4]{3syz}.
On the other hand,  for line arrangements we have a better bound, namely 
$
d_{s} \leq d-2 ,$
see \cite[Corollary 3.5]{Sch}, and for conic-line arrangement we have
$
d_{s} \leq d-1,
$
see \cite[Corollary 1.5]{CMreg}.

In this note we assume that $C$ is a union of $m\geq 4$ smooth curves
$C_j:f_j=0$ of degree $e_j= \deg C_j$, such that $C$ is a normal crossing divisor. Our first aim is to determine explicit expressions for a minimal set of generators for the graded $S$-module $D_0(f)$ in terms of the factors $f_j$, in particular to obtain the corresponding exponents.
Let $\ell(C)$ be the number of irreducible components $C_j$ of $C$ which are lines, e.g. such that $e_j=1$.
Then one has the following result, see \cite[Theorem 2.3, Corollary 2.4 and Corollary 5.2]{Bernd2}.
\begin{thm}
\label{thmA}
With the above notation,  $C$ is an $s$-syzygy curve, where $s$ is given by the following formulas.
\begin{enumerate}

\item If $\ell(C) \leq 2$, then $s= m+2-\ell(C)$.

\item If $\ell(C) \geq 3$, then $s= m-1$.

\end{enumerate}

\end{thm}
We recall now some results from \cite{Edin}. Let $\Omega^j$ be the $S$-module of polynomial differential $j$-forms on $\C^3$ with the usual grading
$$|x|=|y|=|z|=|\dd x|=|\dd y|=|\dd z|=1,$$
for $j=0,1,2,3$. Consider the Koszul complex $(K_f^*,\dd f)$, namely
$$0 \to \Omega^0 \to \Omega^1 \to \Omega^2 \to \Omega^3 \to 0,$$
where the differentials are given by the cup product by the differential $\dd f$ of the polynomial $f$, which is a 1-form of degree $d$, that is
$$\dd f \in \Omega^1_d.$$
We identify a triple 
$$\rho=(a,b,c) \in S^3$$
with the 2-form
$$\omega=a\dd y \wedge \dd z -b\dd x \wedge \dd z+c\dd x \wedge \dd y,$$
and note that
$$\dd f \wedge \omega=af_x+bf_y+cf_z.$$
It follows that we have an isomorphisms of graded $S$-modules
$$\ker \{ \Omega^2 \to \Omega^3 \}=D_0(f)(-2).$$
On the other hand, the module
$$\im \{ \Omega^1 \to \Omega^2 \}$$
is generated by $\omega^x=\dd f \wedge \dd x, \ \omega^y=\dd f \wedge \dd y$ and $\omega^z=\dd f \wedge \dd z$ which correspond to the Koszul type syzygies
$$\kappa^x=(0,f_z,-f_y), \ \kappa^y=(f_z,0,-f_x) \text{ and } \kappa^z=(f_y,-f_x,0).$$
Since $|\kappa^u|=d+1$ for $u=x,y,z$, it follows that the cohomology group
$H^2(K_f^*,\dd f)$ of the Koszul complex satisfies
\begin{equation}
\label{H1}
H^2(K_f^*,\dd f)_j=D_0(f)_{j-2},
\end{equation}
for any $j \leq d$.
Using \cite[Theorem 4.1]{Edin} it follows that $D_0(f)_j=0$ for $j<d-2$ and $D_0(f)_{d-2}=H^2(K_f^*,\dd f)_d$ is $m-1$-dimensional $\C$-vector space, such that the morphism
$$\theta:H^2(K_f^*,\dd f)_d \to H^1(U, \C),  \  \theta (\omega)=\frac{\Delta(\omega)}{f}$$
is an isomorphism of $\C$-vector spaces, where $U= \PP^2 \setminus C$ and $\Delta: \Omega^j \to \Omega^{j-1}$ denotes the contraction with the Euler vector field, see for instance \cite[Lemma 6.1.15]{STH}.
 Consider the differential forms
\begin{equation}
\label{F1}
\omega_j=\dd f \wedge \frac{\dd f_j}{f_j}=\sum_{i=1,m; i\ne j} \frac{f}{f_if_j} \dd f_i \wedge \dd f_j \in \Omega^2
\end{equation}
for $j=1, \ldots, m$ and note that
$$ \theta (\omega_j)=\frac{\Delta(\omega_j)}{f}=d\frac{\dd f_j}{f_j}-e_j\frac{\dd f}{f}.$$
It follows that 
\begin{equation}
\label{F12}
\sum_{j=1}^m\omega_j=0
\end{equation}
and the forms
$\omega _j$ for $j=1,\ldots,m-1$ give a basis of the $\C$-vector space
$H^2(K_f^*,\dd f)_d$. This last claim follows for instance by noting that
$$e_m \theta (\omega_j)-e_j  \theta (\omega_m)=d\left(e_m \frac{\dd f_j}{f_j}-e_j\frac{\dd f_m}{f_m}\right),$$
and using the basis for $H^1(U,\C)$ given in \cite[Equality (3.1)]{Edin}.
With this notation, we can state now our result.
\begin{thm}
\label{thm1}
With the above notation, the graded $S$-module $D_0(f)$ is generated 
 by the syzygies corresponding to a set $G$ of homogeneous differential  forms in $\Omega^2$
as follows.
\begin{enumerate}

\item If $\ell(C)=0$, then $G$ consists of the $m-1$ forms $\omega_j$ for $j=1, \ldots, m-1$ and the three forms $\omega^x$, $\omega^y$ and $\omega^z$. In particular, the exponents of $C:f=0$ in this case are 
$$d_1=\ldots=d_{m-1}=d-2, d_m=d_{m+1}=d_{m+2}=d-1.$$

\item If $\ell(C)=1$, let's say $f_1=x$, then $G$ consists of
 the $m-1$ forms $\omega_j$ for $j=1, \ldots, m-1$ and the two forms $\omega^y$ and $\omega^z$. In particular, the exponents of $C:f=0$ in this case are 
$$d_1=\ldots=d_{m-1}=d-2, d_m=d_{m+1}=d-1.$$

\item If $\ell(C)=2$, let's say $f_1=x$ and $f_2=y$, then $G$ consists of
 the $m-1$ forms $\omega_j$ for $j=1, \ldots, m-1$ and the form  $\omega^z$. In particular, the exponents of $C:f=0$ in this case are 
$$d_1=\ldots=d_{m-1}=d-2, d_m=d-1.$$

\item If $\ell(C) \geq 3$, then $G$ consists of
the $m-1$ forms $\omega_j$ for $j=1, \ldots, m-1$. In particular, the exponents of $C:f=0$ in this case are 
$$d_1=\ldots=d_{m-1}=d-2.$$

\end{enumerate}

\end{thm}

The exponents in Case (4) in the above Theorem can be regarded as a very special case of \cite[Corollary 2.4]{Bernd2} in view of the shift in degrees involved in the isomorphism given by \cite[Proposition 5.1]{Bernd2}.

In fact we can describe the minimal resolution of the graded $S$-module $D_0(f)$ in all the cases covered by Theorem \ref{thm1}. 
The minimal resolution we give below in Theorem \ref{thm10} can be regarded as a very special case of \cite[Theorem 3.6]{BRSS}, see \cite[Example 3.7]{BRSS}. The only plus in our result is the fact that in \cite[Theorem 3.6]{BRSS} there is an assumption of genericity on the forms $f_i$ which might not be easy to check in a given case, while our assumptions ($C$ nodal and the components $C_i:f_=0$ smooth) are rather computable. 
In the next results we set
$$\omega_n'=-\sum_{j=1}^{m-1}\omega_j.$$

\begin{thm}
\label{thm10}
With the above notation, if $\ell(C)=0$, then the minimal resolution of the graded $S$-module $D_0(f)$  has the following shape
$$0 \to \oplus_{j=1}^mS(-d+2-e_j) \to S(-d+2)^{m-1}\oplus S(-d+1)^3 $$
and hence the  minimal resolution of the Jacobian algebra $M(f)$  has the following shape
$$0 \to \oplus_{j=1}^mS(-2d+3-e_j) \to S(-2d+3)^{m-1}\oplus S(-2d+2)^3 \to S(-d+1)^3 \to S.$$
Moreover, the $m$  syzygies $\rho_j$ of $D_0(f)$ correspond to the following relations among the associated 2-forms $\omega_i$ for $i=1, \ldots, m-1$,
$\omega^x$, $\omega^y$ and $\omega^z$
$$\rho_j: \ f_j\omega_j-(f_{j,x}\omega^x+f_{j,y}\omega^x+f_{j,z}\omega^x)=0,$$
for $j=1,\ldots,m,$
where the form $\omega_m$ has to be replaced by $\omega_n'$.

\end{thm}

\begin{thm}
\label{thm11}
With the above notation, if $\ell(C)=1$ and $f_1=x$, then the minimal resolution of the graded $S$-module $D_0(f)$  has the following shape
$$0 \to \oplus_{j=2}^{m}S(-d+2-e_j) \to S(-d+2)^{m-1}\oplus S(-d+1)^2 $$
and hence the  minimal resolution of the Jacobian algebra $M(f)$  has the following shape
$$0 \to \oplus_{j=2}^mS(-2d+3-e_j) \to S(-2d+3)^{m-1}\oplus S(-2d+2)^2 \to S(-d+1)^3 \to S.$$
Moreover, the $(m-1)$  syzygies $\rho_j$ of $D_0(f)$ correspond to the following relations among the associated 2-forms $\omega_i$ for $i=1, \ldots, m-1$, $\omega^y$ and $\omega^z$
$$\rho_j: \ f_j\omega_j-(xf_{j,x}\omega_1+f_{j,y}\omega^x+f_{j,z}\omega^x)=0,$$
for $j=2,\ldots,m,$
where the form $\omega_m$ has to be replaced by $\omega_n'$.

\end{thm}

\begin{thm}
\label{thm12}
With the above notation, if $\ell(C)=2$, and $f_1=x$ and $f_2=y$, then minimal resolution of the graded $S$-module $D_0(f)$  has the following shape
$$0 \to \oplus_{j=3}^mS(-d+2-e_j) \to S(-d+2)^{m-1}\oplus S(-d+1) $$
and hence the  minimal resolution of the Jacobian algebra $M(f)$  has the following shape
$$0 \to \oplus_{j=3}^mS(-2d+3-e_j) \to S(-2d+3)^{m-1}\oplus S(-2d+2) \to S(-d+1)^3 \to S.$$
Moreover, the $(m-2)$  syzygies $\rho_j$ of $D_0(f)$ correspond to the following relations among the associated 2-forms $\omega_i$ for $i=1, \ldots, m-1$ and $\omega^z$
$$\rho_j: \ f_j\omega_j-(xf_{j,x}\omega_1+yf_{j,y}\omega_2+f_{j,z}\omega^x)=0,$$
for $j=3,\ldots,m,$
where the form $\omega_m$ has to be replaced by $\omega_n'$.

\end{thm}

\begin{thm}
\label{thm13}
With the above notation, if $\ell(C)\geq 3$, and $f_1=x$, $f_2=y$ and $f_3=z$, then minimal resolution of the graded $S$-module $D_0(f)$  has the following shape
$$0 \to \oplus_{j=4}^mS(-d+2-e_j) \to S(-d+2)^{m-1} $$
and hence the  minimal resolution of the Jacobian algebra $M(f)$  has the following shape
$$0 \to \oplus_{j=4}^mS(-2d+3-e_j) \to S(-2d+3)^{m-1} \to S(-d+1)^3 \to S.$$
Moreover, the $(m-3)$  syzygies $\rho_j$ of $D_0(f)$ correspond to the following relations among the associated 2-forms $\omega_i$ for $i=1, \ldots, m-1$
$$\rho_j: \ f_j\omega_j-(xf_{j,x}\omega_1+yf_{j,y}\omega_2+zf_{j,z}\omega_3)=0,$$
for $j=4,\ldots,m,$
where the form $\omega_m$ has to be replaced by $\omega_n'$.

\end{thm}
When $m=4$, Theorem \ref{thm13} is a special case of results from \cite{DSres}.

\begin{rk}
\label{rk1}
(i) Since the exponents in Theorem \ref{thm1} have at most 2 values, it is more convenient to denote them as follows:
$$(d-2)_{m-1}(d-1)_3 
\text{ in case (1), }$$
$$(d-2)_{m-1}(d-1)_2
\text{ in case (2), }$$
$$(d-2)_{m-1}(d-1)_1
\text{ in case (3), }$$
and
$$(d-2)_{m-1}
\text{ in case (4) }$$
respectively. \\

(ii) It seems difficult to extend Theorem \ref{thm1} to the case of more than 3 variables, that is to compute the exponents in the $n$-dimensional case. Using the SINGULAR package \cite{Sing}, one may check that the exponents of the surface arrangement $V:f=f_1 \ldots f_m=0$
with
$$f_j= x^p+2^jy^p+3^jz^p+4^jw^p,$$
for $5 \leq m \leq 7$ and $2 \leq p \leq 5$
are, with the notation from the point (i) above,
$$(d-3)_{(m-1)(m-2)/2}(d-2)_{4(m-1)}(d-1)_6.$$
Note that
$$\frac{(m-1)(m-2)}{2}+4(m-1)+6=\binom{m+3}{2},$$
as predicted by \cite[Theorem 3.6]{BRSS} (assuming this arrangement is generic according to the definition given there) or \cite[Corollary 5.2]{Bernd2}. 
It would be nice to have at least a formula for the exponents
in the $n$-dimensional case.
When $V$ is the union of $n+1$ hyperplanes and a smooth hypersurface of degree $>1$ in $\PP^n$, with $V$ a normal crossing divisor, then the coresponding minimal resolution is described in
\cite{DSres}.

(iii) The basis $\omega_j$  of the $\C$-vector space
$H^2(K_f^*,\dd f)_d$ used in this note is different from the basis $\gamma_j$ which occurs in \cite[Theorem 4.1]{Edin}. This new basis has the advantage that it is given by simpler formulas and the secondary syzygies described in Theorem \ref{thm10} have a particularly simple form with respect to it.

(iv)  Theorem \ref{thmA}, Theorem \ref{thm1} and Theorem \ref{thm10}
all fail if the curve $C$ is nodal but not all of its components $C_i$ are smooth.  Using the SINGULAR package \cite{Sing}, one may check that the  curve $C:f=f_1 \ldots f_4=0$
with $f_1=x^3+y^3-3xyz=0$ a nodal cubic and
$$f_j= x^2+2^{j-1}y^2+3^{j-1}z^2=0$$
smooth conics for $j=2,3,4$ is nodal, and its exponents are
$$(7,7,7,8,8,8,8)=(7_3,8_4).$$
In general, the only information we have when the nodal curve $C$ has $m\geq 2$ possibly non-smooth components is that $C$ is an $s$-syzygy curve with $s \geq m-1$ such that
$$d_1=\ldots =d_{m-1}=d-2$$
and the remaining exponents $d_j$, if any, satisfy $d_j>d-2$, see 
\cite[Theorem 4.1]{Edin}. Moreover, by this result, the syzygies 
$$r_1, \ldots, r_{m-1} \in D_0(C)_{d-2}$$
can be chosen to correspond to the 2-forms $\omega_j$ defined in \eqref{F1}. 
\end{rk}

\section{Proof of Theorem \ref{thm1}} 

In this proof, we identify a syzygy $\rho \in D_0(f)$ and the corresponding 2-form in $\Omega^2$.

\medskip
\noindent {\bf Case $\ell(C)=0.$}
Since $\omega_j$ for $j=1, \ldots, m-1$ give a basis for the vector space $D_0(f)_{d-2}$, it follows that they can be chosen as the beginning of a set $G$ of generators for the module $D_0(f)$. To show that
the 3 forms $\omega^x, \ \omega^y, \ \omega^z$ must be added to any such set $G$, it is enough to show that any equality
\begin{equation}
\label{F2}
\sum_{j=1}^{m-1} \ell_j \omega_j+a\omega^x+b\omega^y+c\omega^z=0,
\end{equation}
where $\ell_j \in S_1$ are linear forms and $a,b,c \in \C$, implies
$$\ell_j=0 \text { for all } j=1, \ldots , m-1 \text{ and } a=b=c=0.$$
Since
$$a\omega^x+b\omega^y+c\omega^z=\dd f \wedge \dd \ell,$$
where $\ell=ax+by+cz$, we can replace \eqref{F2} by
\begin{equation}
\label{F3}
\sum_{j=1}^{m-1} \ell_j \omega_j+\dd f \wedge \dd \ell=0.
\end{equation}
We apply the contraction $\Delta$ to the equality \eqref{F3} and get 
\begin{equation}
\label{F4}
\sum_{j=1}^{m-1} \ell_j \left( d\frac{f}{f_j}\dd f_j-e_j\dd f\right) +df \dd \ell-\ell \dd f=0.
\end{equation}
Since $e_m>1$, the only terms not divisible by $f_m$ are 
$$-(\sum_{j=1}^{m-1} e_j\ell_j + \ell) \dd f$$
and hence 
\begin{equation}
\label{F41}
\sum_{j=1}^{m-1} e_j\ell_j + \ell=0.
\end{equation}
On the other hand, the terms not divisible by $f_j$ for $j<m$
are the terms in the expression
$$((d-e_j)\ell_j-\ell)\frac{f}{f_j}\dd f_j,$$
since the partial derivatives of $f_j$ are not divisible by $f_j$, the curve
$C_j:f_j=0$ being smooth. Therefore
\begin{equation}
\label{F42}
(d-e_j)\ell_j-\ell=0.
\end{equation}
It follows using \eqref{F41} and \eqref{F42} that 
\begin{equation}
\label{F40}
\ell_j=0 \text{ for any } j=1, \ldots, m-1  \text{ and } \ell =0.
\end{equation}
This completes the proof of Theorem \ref{thm1} (1).
Indeed, by Theorem \ref{thmA} (1) we know that the total number of syzygies in a minimal generating system is $m+2=(m-1)+3$.

\medskip
\noindent {\bf Case $\ell(C)=1$ and $f_1=x$.}
Note that in this case $\kappa^x= x \kappa^x_0$, where $\kappa^x_0 \in D_0(f)_{d-2}$ is a linear combination of the $\omega_j$'s. Therefore
$\kappa^x$ is not in a minimal set of generators $G$. To complete the proof in this case, we must show that
an equality \eqref{F3} with the linear form $\ell$ not involving $x$ yields
again the vanishings in \eqref{F40}.
Starting from the equation \eqref{F4} we infer as above \eqref{F41} and the equalities
\eqref{F40} if $\ell_1$ is not divisible by $f_1=x$, and in any case
the equalities \eqref{F42} for $j>1$.
Assume now that $\ell_1$ is divisible by $x$.
Then in the sum \eqref{F4}, only the last term is not divisible by $x$,
and hence we get $\ell=0$. This implies $\ell_j=0$ for $j>1$ using the above observation on  \eqref{F42}, and then finally $\ell_1=0$ as well.
Then the proof is completed as in Case $\ell(C)=0$.

\medskip
\noindent {\bf Case $\ell(C)=2$ and $f_1=x$, $f_2=y$.}
Note that in this case $\kappa^x= x \kappa^x_0$ and $\kappa^y= y \kappa^y_0$, where $\kappa^x_0, \kappa^y_0 \in D_0(f)_{d-2}$ are linear combination of the $\omega_j$'s. Therefore
$\kappa^x$ and $\kappa^y$ are not in a minimal set of generators $G$.
To complete the proof in this case, we must show that
an equality \eqref{F3} with the linear form $\ell$ a multiple of $z$ yields
again the vanishings in \eqref{F40}. Starting from the equation \eqref{F4} we infer as above the equalities
\eqref{F42} for $j>2$ only. As in Case (2), the conditions $\ell_1$ divisible by $x$ or $\ell_2$ divisible by $y$ yield a contradiction, e. g. by using \eqref{F41}.
If we assume that $\ell_1$ is not divisible by $x$ and $\ell_2$ is not divisible by $y$, we get as above that the equalities in \eqref{F42}
hold for $j=1$ and $j=2$ as well.
Then the proof is completed as in Case $\ell(C)=0$.

\medskip
\noindent {\bf Case $\ell(C)\geq 3$}
By Theorem \ref{thmA} (2) we know that a minimal set of generators $G$ has $m-1$ elements. On the other hand we know that $\omega_j \in G$ by our discussion before Theorem \ref{thm1}. Therefore there are no other elements in the generating set $G$ in this case.

\section{Proof of Theorems \ref{thm10}, \ref{thm11},\ref{thm12} and \ref{thm13}} 

\subsection{Proof of Theorems \ref{thm10}}

\noindent{\bf Step 1.} Independence of the syzygies $\rho_j$'s.

We first describe in a new way the $m$  syzygies of $D_0(f)$ described in Theorem \ref{thm10}. 
For $j<m$, the image of $\rho_j$ in $S(-d+2)^{m-1}\oplus S(-d+1)^3 $
is given by 
$$(0, \ldots, 0, f_j,0, \ldots, 0, -f_{j,x},-f_{j,y},-f_{j,z})$$
where $f_j$ is on the $j$-th coordinate.
Moreover, 
the image of $\rho_m$ in $S(-d+2)^{m-1}\oplus S(-d+1)^3 $
is given by 
$$(-f_m, \ldots , -f_m, -f_{m,x},-f_{m,y},-f_{m,z}).$$
Consider the matrix $\M$ with $m$ rows and $m+2$ columns, obtained by taking as $j$-th row the $m+2$ components of the image of $\rho_j$
described above. We claim that 
\begin{equation}
\label{R}
\rank \M=m,
\end{equation}
 in other words
there no non-trivial linear combinations
$$A_1 \rho_1+ \ldots + A_m \rho_m=0,$$
where $A_j \in S$ are (homogeneous) polynomials. 
Indeed, assume that such a linear combination exists and choose one such combination where $\deg A_j+\deg \rho_j$ is minimal. If $A_m=0$, we get a contradiction, since the first $m-1$ rows and the first $m-1$ colons yields a square matrix with determinant $f_1 \ldots f_{m-1} \ne 0$.
Hence we may suppose $A_m \ne 0$
and let $a\geq 0$ be the largest integer such that $f_m^ a$ divides $A_m$. Note that looking at the first $m-1$ columns we get that
$A_j$ is divisible by $f_m^{a+1}$ for any $j<m$. If we look now at the
$m$-th column, we see that $f_{m,x}$ is not divisible by $f_m$,
and hence $A_m$ has to be divisible by $f_m^{a+1}$, a contradiction.
Indeed, $A_mf_{m,x}$ is a sum of terms, each one divisible by $f_m^{a+1}$.

\noindent{\bf Step 2.} Recall of a property of the degrees of second order syzygies for $M(f)$, which are the first order syzygies of $D_0(f)$.

Our curve being nodal of degree $d\geq 6$ is not free, hence
 the minimal resolution for the Milnor algebra $M(f)$ has the form
\begin{equation}
\label{res2A}
0 \to \oplus_{i=1} ^{m}S(-c_i) \to \oplus_{j=1} ^{m+2}S(1-d-d_j)\to S^3(1-d)  \to S,
\end{equation}
with $c_1\leq ...\leq c_{m}$ and $d_1= ...= d_{m-1}=d-2$ and $d_m=d_{m+1}=d_{m=2}=d-1$.
It follows from \cite[Lemma 1.1]{HS} that one has
\begin{equation}
\label{res2B}
c_j=d+d_{j+2}-1+\epsilon_j,
\end{equation}
for $j=1,...,m$ and some integers $\epsilon_j\geq 1$. Using \cite[Formula (13)]{HS}, it follows that one has
\begin{equation}
\label{res2C}
d_1+d_2=d-1+\sum_{i=1} ^{m}\epsilon_j.
\end{equation}
With this notation, we call $\epsilon_j$ the degree of the $j$-th secondary syzygy and note that in our case
$$\sum_{i=1} ^{m}\epsilon_j=d_1+d_2-(d-1)=2(d-2)-(d-1)=d-3.$$
If our claim in Theorem \ref{thm10} holds, then
our second syzygy $\rho_j$ has a degree $\epsilon'_j=e_j$ for
$j<m-2$ and $\epsilon'_j=e_j-1$ for $j=m-2,m-1,m$. It follows that
$$\sum_{i=1} ^{m}\epsilon'_j=d-3.$$

\noindent{\bf Step 3.} The rest of the proof.

Let $r_1, \ldots, r_m$ be a minimal set of homogeneous generators for the  syzygy module of $D_0(f)$ such that $\deg r_i \leq \deg r_j$ for all $i<j$. Since the submodule generated by
$\rho_1, \ldots,\rho_j$ has rank $j$ for any $j=1,\ldots,m$ by Step 1 above, it follows that, up-to a possible reordering of the $r_k$'s, for any $j$, one has
$$\rho_j \notin S[r_1,\ldots, r_{j-1}] \text{ and } \rho_j \in S[r_1,\ldots, r_{j'}],$$
where $S[u,v,w,...]$ denotes the $S$-submodule generated by $u,v,w,...$ and with $j'$ minimal with this property. In particular $j' \geq j$.
These relations imply that
$$\deg r_j \leq \deg r_j' \deg \rho_j$$
for all $j=1,\ldots,m$. Taking the sum over all $j$'s, and using Step 2, we see that in fact one has the equalities
$$\deg r_j=\deg \rho_j$$
for all $j$. It follows that any syzygy $r_j$ is in the submodule $S[\rho_1,\ldots,\rho_j$, and therefore
$\rho_1, \ldots,\rho_m$ is also a minimal set of generators for the  syzygy module of $D_0(f)$.

\subsection{Proof of Theorems \ref{thm11}} We indicate only the changes to be done in the above proof to treat this case.
Note that the first  syzygy $\rho_1$ in Theorem \ref{thm10} becomes
$$x\omega_1=\omega^x,$$
and hence $\omega^x$ has to be discarded from the minimal set of generators of $D_0(f)$ and $\rho_1$ from the set of minimal generators of the  syzygy module of $D_0(f)$. Consider the new matrix $\M^x$ obtained from the matrix $\M$ by the following operations. First, replace the first column $C^1$ by the difference $C^1-xC^{m}$, where
$C^{m}$ denotes the $m$-th column in $\M$. Then discard the first row (which has now only one non-zero entry) and the $m$-th column in $\M$ to get the matrix $\M^x$. It is clear that
\begin{equation}
\label{Rx}
\rank \M^x=\rank \M-1=m-1.
\end{equation}
Then the proof is completed as in the case of Theorem \ref{thm10}.

\subsection{Proof of Theorems \ref{thm12}} 

The first two syzygies $\rho_1$ and $\rho_2$ in Theorem \ref{thm10} become
$$x\omega_1=\omega^x \text{ and } y \omega _2=\omega^y,$$
and hence $\omega^x$ and $\omega^y$ have to be discarded from the minimal set of generators of $D_0(f)$ and $\rho_1$ and $\rho_2$ from the set of minimal generators of the syzygy module of $D_0(f)$.
Consider the new matrix $\M^{xy}$ obtained from the matrix $\M^x$ by the following operations. First, replace the first column $C^{x1}$ by the difference $C^{x1}-yC^{xm}$, where
$C^{xm}$ denotes the $m$-th column in $\M^x$. Then discard the first row (which has now only one non-zero entry) and the $m$-th column in $\M^x$ to get the matrix $\M^{xy}$. It is clear that
\begin{equation}
\label{Rxy}
\rank \M^{xy}=\rank \M^x-1=m-2.
\end{equation}
Then the proof is completed as in the case of Theorem \ref{thm10}.

\subsection{Proof of Theorems \ref{thm13}} 

The first three  syzygies $\rho_1$, $\rho_2$ and $\rho_3$ in Theorem \ref{thm10} become
$$x\omega_1=\omega^x, \   y\omega_2=\omega^y   \text{ and } z \omega _3=\omega^z,$$
and hence $\omega^x$,  $\omega^y$ and $\omega^z$ have to be discarded from the minimal set of generators of $D_0(f)$ and $\rho_1$, $\rho_2$ and $\rho_3$ from the set of minimal generators of the  syzygy module of $D_0(f)$.
Consider the new matrix $\M^{xyz}$ obtained from the matrix $\M^{xy}$ by the following operations. First, replace the first column $C^{xy1}$ by the difference $C^{xy1}-zC^{xym}$, where
$C^{xym}$ denotes the $m$-th column, which is in fact the last column, in $\M^{xy}$. Then discard the first row (which has now only one non-zero entry) and the $m$-th column in $\M^{xy}$ to get the matrix $\M^{xyz}$. It is clear that
\begin{equation}
\label{Rxyz}
\rank \M^{xyz}=\rank \M^{xy}-1=m-3.
\end{equation}
Then the proof is completed as in the case of Theorem \ref{thm10}.

\end{document}